\pdfoutput=1
\RequirePackage{ifpdf}
\ifpdf 
\documentclass[pdftex]{sigma}
\else
\documentclass{sigma}
\fi

\usepackage{mathrsfs}

\def\A{{\mathcal A}}

\def\th{\theta}
\def\Th{\Theta}

\def\N{\mathbb{N}}

\def\X{\mathscr X}

\def\QQ{\mathfrak Q}

\def\AA{\mathfrak A}

\def\C{\mathcal{C}}

\def\B{{\mathcal B}}

\def\K{\mathcal K}

\def\Q{\mathcal Q}
\def\P{\mathcal P}
\def\R{\mathcal{R}}

\def\Ai{\mathcal A^\infty}
\def\AAi{\mathfrak A^\infty}

\def\BB{\mathfrak B}
\def\CC{{\mathfrak C}}
\def\PP{{\mathfrak P}}
\def\RR{{\mathfrak R}}

\def\J{\mathcal J}
\def\JJ{\mathfrak J}

\def\I{{\rm 1\kern-.26em I}}

\def\si{\sigma}
\def\Si{\Sigma}

\def\1{\mathfrak{1}}
\def\0{\mathfrak{0}}

 \def\hb{\hbar}

\def\<{\langle}
\def\>{\rangle}

\providecommand{\CC}{\mathfrak{C}}

\newtheorem{Theorem}{Theorem}[section]
\newtheorem{Lemma}[Theorem]{Lemma}

\newtheorem{Corollary}[Theorem]{Corollary}
\newtheorem{Proposition}[Theorem]{Proposition}
{\theoremstyle{definition}
\newtheorem{Definition}[Theorem]{Definition}
\newtheorem{Example}[Theorem]{Example}
\newtheorem{Remark}[Theorem]{Remark}
}

\begin{document}

\allowdisplaybreaks

\renewcommand{\PaperNumber}{091}

\FirstPageHeading

\ShortArticleName{Covariant Fields of $C^*$-Algebras under Rief\/fel Deformation}

\ArticleName{Covariant Fields of $\boldsymbol{C^*}$-Algebras\\ under Rief\/fel Deformation}

\Author{Fabian BELMONTE~$^\dag$ and Marius M\u{A}NTOIU~$^\ddag$}

\AuthorNameForHeading{F.~Belmonte and M.~M\u{a}ntoiu}

\Address{$^\dag$~728, building A, SISSA-ISAS, Via Bonomea 265, 34136 Trieste, Italy}
\EmailD{\href{mailto:fabian.belmonte@sissa.it}{fabian.belmonte@sissa.it}}

\Address{$^\ddag$~Departamento de Matem\'aticas, Universidad de Chile,\\
\hphantom{$^\ddag$}~Las Palmeras 3425, Casilla 653, Santiago, Chile}
\EmailD{\href{mailto:mantoiu@uchile.cl}{mantoiu@uchile.cl}}

\ArticleDates{Received August 26, 2012, in f\/inal form November 22, 2012; Published online November 28, 2012}

\Abstract{We show that Rief\/fel's deformation sends covariant $\C(T)$-algebras into $\C(T)$-algebras.
We also treat the lower semi-continuity issue, proving that Rief\/fel's deformation transforms
covariant continuous f\/ields of $C^*$-algebras into continuous f\/ields of $C^*$-algebras.
Some examples are indicated, including certain quantum groups.}

\Keywords{pseudodif\/ferential operator; Rief\/fel deformation; $C^*$-algebra; continuous f\/ield;
noncommutative dynamical system}

\Classification{35S05; 81Q10; 46L55; 47C15}

\section{Introduction}\label{duci}

Let $T$ be a locally compact topological space, always assumed to be Hausdorf\/f. We denote by~$\C(T)$ the Abelian $C^*$-algebra
of  all complex continuous functions on~$T$ that decay at inf\/inity (are arbitrarily small outside large compact subsets).
A $\C(T)$-algebra~\cite{Bl,Ni,Wi} is a $C^*$-algebra~$\B$ together with a nondegenerated injective morphism from $\C(T)$ to
the center of~$\B$ (multipliers are used if~$\B$ is not unital).
The main role of the
concept of $\C(T)$-algebra consists in translating in a simple and ef\/f\/icient way the idea that $\B$ is f\/ibred in the sense of
$C^*$-algebras over the base~$T$~\cite{Fe,To1}.
Actually $\C(T)$-algebras can be seen as upper-semi-continuous f\/ields of $C^*$-algebras over the base~$T$.
Lower-semi-continuity can also be put in this setting if one also uses the space of all primitive ideals \cite{Lee,Ni,RW,Rie2,Wi};
we refer to~\cite{Bl2} for deeper results involving semi-continuity.
We intend to put these concepts in the perspective of Rief\/fel quantization.

Rief\/fel's calculus \cite{Rie1,Rie2} is a machine that applies to $C^*$-dynamical systems and their morphisms.
The necessary ingredients are an action of the vector group $\Xi:=\mathbb R^d$ by automorphisms of a $C^*$-algebra as well
as a skew symmetric linear operator of~$\Xi$.
When morphisms are involved, they are always assumed to intertwine the existing actions.

Rief\/fel's machine is actually meant to be a quantization.
The initial data def\/ine a natural Poisson structure,
regarded as a mathematical modelization of the observables of a classical physical system.
After applying the machine to this classical data one gets a deformed $C^*$-algebra seen as the family of observables of the same system,
but written in the language of quantum mechanics. By varying a convenient parameter (Planck's constant $\hbar$)
one can recover the Poisson structure (at $\hbar=0$) from the $C^*$-algebras def\/ined at $\hbar\ne0$ in a way that satisf\/ies
certain natural axioms~\cite{La,Rie1,Rie2}.

The spirit of this deformation procedure is that of a pseudodif\/ferential
theory~\cite{Fol89}. At least in simple situations the multiplication in the initial $C^*$-algebra is just point-wise multiplication
of functions def\/ined on some locally compact topological space, on which~$\Xi$ acts by homeomorphisms.
The noncommutative product in the quantized algebra can be interpreted as a symbol composition of a pseudodif\/ferential type.
Actually the concrete formulae generalize and are motivated by the usual Weyl calculus.

In a setting where all the relevant concepts make sense, we prove in Theorem~\ref{crifel} and Proposition~\ref{pieligroasa} their compatibility: by Rief\/fel quantization an upper-semi-continuous f\/ields of $C^*$-algebras
is turned into an upper-semi-continuous f\/ields of $C^*$-algebras with f\/ibers which are easy to identify; the proof
uses $\C(T)$-algebras. Finally, using primitive ideals techniques,
we show the analog of this result for lower-semi-continuity; the key technical result is Proposition~\ref{segund}.
Putting everything together one gets

\begin{Theorem}\label{mein}
Rieffel quantization transforms covariant continuous fields of $C^*$-algebras into
covariant continuous fields of $C^*$-algebras.
\end{Theorem}

We illustrate the result by some examples in Section~\ref{seretide}. Most of them involve an
Abelian initial algebra $\A$. In this case the information is encoded in a topological dynamical system
with locally compact space $\Si$ and the upper-semi-continuous f\/ield property can be read in the existence of a continuous
covariant surjection $q:\Si\rightarrow T$; if this one is open, then lower-semi-continuity also holds.
If the orbit space of the dynamical system is Hausdorf\/f, it serves as a good space $T$ over which the Rief\/fel
deformed algebra can be decomposed, with easily identif\/ied f\/ibers. This can be
used to show that the $C^*$-algebras of some
compact quantum groups constructed in~\cite{Rie3} can be written as continuous f\/ields, some of the f\/ibers being isomorphic
to certain noncommutative tori.

Eventual connections of the present work to the model of quantum spacetime presented in~\cite{DFR} and~\cite{Pi} will be investigated elsewhere, hopefully.

After the present work was completed, we learnt of the recent articles \cite{HM1,HM2} in which a statement comparable with our
Theorem~\ref{crifel} is included (essentially without proof) with interesting applications. These papers reply on Kasprzak's
reformulation~\cite{Ka} of the Rief\/fel deformation (based on Landstad's characterization of crossed products). However, in \cite{HM1,HM2} there is no comment about how nondegenerecy is proved, and this is crucial in the def\/inition of a $\mathcal C(T)$-algebra (cf.\ Def\/inition~\ref{cedete}). It is easy to give counterexamples showing how important this condition is for the theory. Moreover, for us this was the main technical problem in proving Theorem~\ref{crifel}; we solved it using the highly nontrivial Dixmier--Malliavin theorem~\cite{DM}. Our detailed proof
based on the more traditional approach is intended to support future work on spectral theory of pseudo\-dif\/ferential operators
(see also \cite{BM,Ma}), as well as applications to quantum mechanics.
The lower-semi-continuity part (contained in Section~\ref{isra}), often described as dif\/f\/icult, is not treated in \cite{HM1,HM2} and seems not to have
correspondence in the literature.

Let us mention that a dif\/ferent proof of the main result (Theorem~\ref{mein}) has been given in~\cite{BM}, where this result is needed as an important ingredient in proving spectral continuity properties of quantum Hamiltonians def\/ined as phase-space anisotropic pseudodif\/ferential operators. That paper has been submitted after the completion of the present article and refers to it explicitly. In~\cite{BM} the proof relies on a recently found connection between Rief\/fel deformation and twisted crossed products~\cite{BMa} and on the deep analysis from~\cite{Rie}, therefore being less direct or self-contained than the present one.

\section{Rief\/fel's pseudodif\/ferential calculus: a short review}\label{sectra}

We start by describing brief\/ly Rief\/fel deformation \cite{Rie1,Rie2} in a slightly restricted setting.
The initial object, containing {\it the classical data}, is a quadruplet
$\left(\A,\Theta,\Xi,[\![\cdot,\cdot,]\!]\right)$. The pair $\left(\Xi,[\![\cdot,\cdot]\!]\right)$ will
be taken to be a $2n$-dimensional symplectic vector space.
On the other hand $\left(\A,\Theta,\Xi\right)$ is a $C^*$-{\it dynamical system},
meaning that the vector group
acts strongly continuously by automorphisms of the (possibly noncommutative) $C^*$-algebra $\A$.
Let us denote by $\Ai$ the family of elements $f$ such that the mapping
$\Xi\ni X\mapsto\Th_X(f)\in\mathcal A$ is
$C^\infty$. It is a dense $^*$-algebra of $\A$ and also a Fr\'echet algebra with the family of semi-norms
\begin{gather}\label{semicar}
\| f\|_\A^{(k)}:=\sum_{|\alpha|\le
k}\frac{1}{|\alpha|!}\|\partial_X^\alpha\left[\Th_X(f)\right]_{X=0}\|_\A \equiv\sum_{|\alpha|\le
k}\frac{1}{|\alpha|!}\|\delta^\alpha (f)\|_\A ,
\qquad k\in\N .
\end{gather}
To quantize the above structure, one keeps the involution unchanged but introduces on $\Ai$ the product
\begin{gather*}
f\,\#\,g:=\pi^{-2n}\int_\Xi\int_\Xi dYdZ\,e^{2i[\![Y,Z]\!]} \Th_Y(f) \Th_Z(g) ,
\end{gather*}
suitably def\/ined by oscillatory integral techniques~\cite{Fol89,Rie1}. One gets a $^*$-algebra $(\Ai,\#,^*)$, which admits a
$C^*$-completion~$\AA$ in a $C^*$-norm $\|\cdot\|_\AA$ def\/ined by Hilbert module techniques~\cite{Rie1}.
The action~$\Th$ leaves~$\Ai$ invariant and extends to a strongly continuous action of the
$C^*$-algebra~$\AA$, that will also be denoted by~$\Th$. The space $\AAi$ of $C^\infty$-vectors coincides with~$\Ai$
and it is a~Fr\'echet space with semi-norms
\begin{gather}\label{semon}
\| f\|_\AA^{(k)}:=\sum_{|\alpha|\le k}\frac{1}{|\alpha|!}\|\partial_X^\alpha\left[\Th_X(f)\right]_{X=0}
\|_\AA \equiv\sum_{|\alpha|\le
k}\frac{1}{|\alpha|!}\|\delta^\alpha (f)\|_\AA ,\qquad k\in\N  .
\end{gather}
By Proposition 4.10 in \cite{Rie1}, there exist $k\in\N $ and $C_k>0$ such that
$\| f\|_\AA \le C_k\| f\|^{(k)}_\A$ for any $f\in\A^\infty=\AA^\infty$.
Replacing here $f$ by $\delta^\alpha f$ for every multi-index $\alpha$, it follows that on $\A^\infty$ the topology given by
the semi-norms~(\ref{semicar}) is f\/iner than the one given by the semi-norms~(\ref{semon}). As a~consequence of
Theorem~7.5 in~\cite{Rie1}, the role of the $C^*$-algebras~$\A$ and $\AA$ can be reversed: one obtains~$\A$ as the deformation
of~$\AA$ by replacing the skew-symmetric form $[\![\cdot,\cdot]\!]$ by $-[\![\cdot,\cdot]\!]$.
Thus~$\A^\infty$ {\it and $\AA^\infty$ coincide as Fr\'echet spaces}.

The deformation transfers to $\Xi$-morphisms. Let
$\left(\A_j,\Theta_j,\Xi,[\![\cdot,\cdot]\!]\right)$, $j=1,2$, be two classical data and
let $\R:\A_1\to\A_2$ be a $\Xi$-morphism, i.e.\ a $C^*$-morphism intertwining the
two actions~$\Theta_1$,~$\Theta_2$. Then $\R$ sends $\A^\infty_1$ into $\A^\infty_2$ and extends to a morphism
$\RR:\AA_1\to\AA_2$ that also intertwines the corresponding actions. In this way, one obtains a
covariant functor.
The functor is exact~\cite{Rie1}: it preserves short exact sequences of~$\Xi$-morphisms. Namely, if~$\J$ is a~(closed,
self-adjoint, two-sided) ideal in $\A$ that is invariant under $\Theta$, then its deformation $\JJ$ can be
identif\/ied with an invariant ideal in~$\AA$ and the quotient~$\AA/\JJ$ is canonically isomorphic to the
deformation  of the quotient~$\A/\J$  under the natural quotient action.

We will refer to {\it the Abelian case} under the following circumstances:
A continuous action $\Theta$ of $\Xi$ by homeomorphisms of the locally compact Hausdorf\/f space~$\Sigma$ is given. For
$(\sigma,X)\in\Sigma\times\Xi$ we are going to use all the notations
$\Theta(\sigma,X)=\Theta_X(\sigma)=\Theta_\sigma(X)\in\Sigma$ for the $X$-transformed of the point~$\sigma$.
The function $\Theta$ is continuous and the homeomorphisms
$\Th_X$, $\Th_Y$ satisfy $\Th_X\circ\Th_Y=\Th_{X+Y}$ for every $X,Y\in\Xi$.

We denote by $\C(\Sigma)$ the Abelian $C^*$-algebra of all complex continuous functions on $\Sigma$
that are arbitrarily small outside large compact subsets of $\Sigma$. When~$\Sigma$ is compact,
$\C(\Si)$ is unital. The action $\Theta$ of $\Xi$ on $\Sigma$ induces an action of~$\Xi$ on
$\C(\Si)$ (also denoted by~$\Theta$) given by~$\Theta_X(f):=f\circ\Theta_X$.
This action is strongly continuous, i.e.\ for any $f\in\C(\Sigma)$ the mapping
\begin{gather}\label{ciuca}
\Xi\ni X\mapsto\Theta_X(f)\in\C(\Sigma)
\end{gather}
is continuous; thus we are placed in the setting presented above.
We denote by $\C(\Si)^\infty\equiv\C^\infty(\Sigma)$ the set of elements $f\in\C(\Sigma)$ such that the
mapping \eqref{ciuca} is $C^\infty$; it is a dense $^*$-algebra of~$\C(\Sigma)$.
The general theory supplies a $C^*$-algebra $\AA\equiv\CC(\Si)$ (noncommutative in general),
acted continuously by the group~$\Xi$, with
smooth vectors $\CC^\infty(\Si)=\C^\infty(\Si)$.

\section[Families of $C^*$-algebras]{Families of $\boldsymbol{C^*}$-algebras}\label{secintra}

Now we give a short review of $\C(T)$-algebras and semi-continuous f\/ields of $C^*$-algebras
(see \cite{Bl,La,Lee,Ni,Rie,Wi} and references therein), outlining the connection between the two notions.

If $\B$ is a $C^*$-algebra, we denote by $\mathcal M(\B)$ its multiplier algebra and by $\mathcal Z\mathcal M(\B)$ its center.
If $\mathcal B_1$, $\mathcal B_2$ are two vector subspaces of $\mathcal M(\B)$, we set $\B_1\cdot\B_2$ for the vector
subspace generated by $\{b_1 b_2\mid b_1\in\B_1,b_2\in\B_2\}$.
We are going to denote by $\C(T)$ the $C^*$-algebra of all complex continuous functions on the (Hausdorf\/f)
locally compact space $T$ that decay at inf\/inity.

\begin{Definition}\label{cedete}
We say that $\B$ is a {\it $\C(T)$-algebra}
if a nondegenerate monomorphism $\Q:\C(T)\rightarrow\mathcal Z\mathcal M(\B)$ is given.
\end{Definition}

We recall that nondegeneracy means that the ideal $\Q[\C(T)]\cdot\B$
is dense in~$\B$.

\begin{Definition}\label{niulaif}
By {\it upper-semi-continuous field of $C^*$-algebras} we mean a family of epimorphisms of $C^*$-algebras
$\big\{\mathcal B\overset{\P(t)}{\longrightarrow}\mathcal B(t)\mid t\in T\big\}$
indexed by a locally compact topological space~$T$ and satisfying:
\begin{enumerate}\itemsep=0pt
\item
For every $b\in\B$ one has $\| b\|_\B=\sup_{t\in T}\| \P(t)b\|_{\B(t)}$.
\item
For every $b\in\B$ the map
$T\ni t\mapsto\| \P(t)b\|_{\B(t)}$ is upper-semi-continuous and decays at inf\/inity.
\item
There is a multiplication $\C(T)\times\B\ni(\varphi,b)\rightarrow\varphi\ast b\in\B$ such that
\[
\P(t)[\varphi\ast b]=\varphi(t) \P(t)b ,\qquad \forall\,t\in T,\quad \varphi\in\C(T) ,\quad  b\in\B .
\]
\end{enumerate}
If, in addition, the map $t\mapsto\| \P(t)b\| $ is continuous for every $b\in \B$, we say that
\[
\big\{\mathcal B\overset{\P(t)}{\longrightarrow}\mathcal B(t)\mid t\in T\big\}
\] is {\it a~continuous field of
$C^*$-algebras}.
\end{Definition}

The requirement~2 is clearly equivalent with the condition that
for every $b\in\B$ and every $\epsilon>0$ the subset $\{t\in T\mid\| \P(t)b\|_{\B(t)}\ge \epsilon\}$  is compact.
One can rephrase~1 as $\cap_t\ker[\P(t)]=\{0\}$,
so one can identify $\B$ with a $C^*$-algebra of sections of the f\/ield; this make the connection with other approaches,
as that of~\cite{Ni} for example. It will always be assumed that $\B(t)\ne\{0\}$ for all $t\in T$.

We are going to describe brief\/ly in which way the two def\/initions above are actually equivalent.

First let us assume that $\B$ is a $\C(T)$-algebra and denote by $\C_t(T)$ the ideal of all the functions in
$\C(T)$ vanishing at the point $t\in T$. We get ideals
$\mathcal I(t):=\overline{\Q\left[\C_t(T)\right]\cdot\B}$ in $\B$, quotients $\B(t):=\B/\mathcal I(t)$ as well as canonical
epimorphisms $\P(t):\B\rightarrow \B(t)$. One also sets
\begin{gather*}
\varphi\ast b:=\Q(\varphi)b ,\qquad \forall\,\varphi\in\C(T),\quad b\in\B .
\end{gather*}
Then $\big\{\mathcal B\overset{\P(t)}{\longrightarrow}\mathcal B(t)\mid t\in T\big\}$ is an upper-semi-continuous
f\/ield of $C^*$-algebras with multiplication~$\ast$.

Conversely, if an upper-semi-continuous f\/ield
$\big\{\mathcal B\overset{\P(t)}{\longrightarrow}\mathcal B(t)\mid t\in T\big\}$ is given,
also involving the multiplication~$\ast$, we set
\begin{gather*}
\Q: \ \C(T)\rightarrow\mathcal Z\mathcal M(\B),\qquad \Q(\varphi)b:=\varphi\ast b.
\end{gather*}
In this way one gets a $\C(T)$-algebra and each of the quotients $\B/\mathcal I(t)$ is isomorphic to the f\/i\-ber~$\B(t)$.

To discuss lower semi-continuity we need ${\rm Prim}(\B)$, the space of all the primitive ideals
(kernels of irreducible representations) of $\B$. The hull-kernel topology
turns ${\rm Prim}(\B)$ into a locally compact (not necessarily Hausdorf\/f) topological space. We recall that {\it the hull application}
$\mathcal J\mapsto h(\mathcal J):=\{\mathcal K\in {\rm Prim}(\B)\mid \mathcal J\subset\mathcal K\}$
realizes a containment reversing bijection between the family of ideals of $\B$ and the family
of closed subsets of ${\rm Prim}(\B)$. Its inverse is {\it the kernel map}
$\Omega\mapsto k(\Omega):=\cap_{\mathcal K\in\Omega}\mathcal K$, which is also decreasing.

The Dauns--Hof\/fman theorem establishes the existence of a unique isomorphism
\[
\Gamma: \ BC[{\rm Prim}(\B)]\to \mathcal Z\mathcal M(\B) ,
\] where $BC[{\rm Prim}(\B)]$ is the $C^*$-algebra of bounded and continuous
functions over ${\rm Prim}(\B)$, such that for each $\K\in {\rm Prim}(\B)$,
$\Psi\in BC[{\rm Prim}(\B)]$ and $b\in\B$ we have $\Gamma(\Psi)b+\mathcal K=\Psi(\mathcal K)b+\mathcal K$.
For a~detailed study of the space ${\rm Prim}(\B)$ and a
proof of the Dauns--Hof\/fman theorem, cf.\ Sections~A.2 and~A.3 in~\cite{RW}.
Let us suppose that there is a continuous surjective function $q:{\rm Prim}(\B)\to T$. Then we can def\/ine
$\Q:\C(T)\to \mathcal Z\mathcal M(\B)$ by $\Q(\varphi)=\Gamma(\varphi\circ q)$ and one can check that
$\Q$ endows $\B$ with the structure of a $\C(T)$-algebra.

On the other hand, if we have a nondegenerate monomorphism $\Q:\C(T)\to \mathcal Z\mathcal M(\B)$, we can def\/ine canonically
a continuous map $q:{\rm Prim}(\B)\to T$. One has $q(\mathcal K)=t$ if and only if $\mathcal I(t)\subset \mathcal K$,
and we can recover $\Q$ from the above construction.
Moreover {\it the map $T\ni t\to\| b(t)\|_{\B(t)}\in\mathbb R_+$ is continuous
for every $b\in\B$ $($so we have a continuous field of $C^*$-algebras$)$ if and only if $q$ is open}.
For the proof of this facts see Propositions~C.5 and~C.10 in~\cite{Wi}.

\section[Covariant $\C(T)$-algebras and upper-semi-continuity under Rieffel quantization]{Covariant $\boldsymbol{\C(T)}$-algebras and upper-semi-continuity\\ under Rief\/fel quantization}\label{isra}

Let $T$ be a locally compact Hausdorf\/f space and $(\A,\Th,\Xi,[\![\cdot,\cdot]\!])$ a classical data.
The canonical $C^*$-dynamical system def\/ined by Rief\/fel quantization is $(\AA,\Th,\Xi)$.

\begin{Definition}\label{cedete+}
We say that $\A$ is a {\it covariant $\C(T)$-algebra with respect to the action} $\Th$
if a~nondegenerate monomorphism $\Q:\C(T)\rightarrow\mathcal Z\mathcal M(\A)$ is given
(so it is a $\C(T)$-algebra) and in addition one has
\begin{gather}\label{incep}
\Th_X[\Q(\varphi)f]=\Q(\varphi)\left[\Th_X(f)\right],\qquad \forall\,f\in\A,\quad X\in\Xi,\quad \varphi\in\C(T) .
\end{gather}
\end{Definition}

We intend to prove that the Rief\/fel quantization transforms covariant $\C(T)$-algebras into covariant $\C(T)$-algebras.
For this and for a further result identifying the emerging quotient algebras, we are going to need

\begin{Lemma}\label{densly}
Let $I$ be an ideal of $\C(T)$ and denote by $\overline{\Q(I)\cdot\A}$ the closure of $\Q(I)\cdot\A$ in the $C^*$-algebra~$\A$.
Then $\Q(I)\cdot\A^\infty$ is dense in
$\big(\overline{\Q(I)\cdot\A}\big)^\infty\equiv\big(\overline{\Q(I)\cdot\A}\big)\cap\A^\infty$
for the Fr\'echet topology inherited from~$\A^\infty$.
\end{Lemma}

\begin{proof}
By the covariance condition $\overline{\Q(I)\cdot\A}$ is an invariant ideal of $\A$.

The proof uses regularization. Consider {\it the integrated form of} $\Th$, i.e.\ for each $\Phi\in C_c^\infty(\Xi)$
(compactly supported smooth function) and $g\in\A$ def\/ine
\[
\Th_\Phi(g)=\int_\Xi dY\Phi(Y)\Th_Y(g).
\]
Note that for every $X\in\Xi$ one has
\[
\Th_X\left[\Th_\Phi(g)\right]=\int_\Xi dY\Phi(Y-X)\Th_Y(g).
\]
Then $\Th_\Phi(g)\in\A^\infty$ and for each multi-index $\mu$ we have
\[
\delta^{\mu}\left[\Th_\Phi(g)\right]=(-1)^{|\mu|}\Th_{\partial^{\mu}\Phi}(g)\qquad \text{and}\qquad
\|\delta^{\mu}\left[\Th_\Phi(g)\right]\|_\A \leq \|\partial^{\mu}\Phi\|_{L^1(\Xi)}
\| g\|_\A .
\]
One of the deepest theorems about smooth algebras, the Dixmier--Malliavin theorem~\cite{DM}, say that~$\A^\infty$
is generated (algebraically) by the set of all the elements of the form  $\Th_\Phi(g)$ with $\Phi\in C_c^\infty(\Xi)$ and
$g\in\A$. Replacing~$\A$ with $\overline{\Q(I)\cdot\A}$,
for $f\in\big(\overline{\Q(I)\cdot\A}\big)^\infty$ there exist $\Phi_1,\dots,\Phi_m\in C_c^\infty(\Xi)$ and
$f_1,\dots,f_m\in\overline{\Q(I)\cdot\A}$ such that $f=\sum\limits_{i=1}^m\Th_{\Phi_i}(f_i)$.
Let $\epsilon >0$ and f\/ix a multi-index $\alpha$. Choose $g_1,\dots ,g_m\in\Q(I)\cdot\A$ such that for each $i$
\[
\| f_i-g_i\|_\A \le\frac{\epsilon}{m\|\partial^\alpha\Phi_i\|_{L^1(\Xi)}} .
\]
Then
\begin{gather*}
\left\| \delta^\alpha \left(f-\sum_{i=1}^m\Th_{\Phi_i}(g_i)\right)\right\|_\A=
\left\|\sum_{i=1}^m\Th_{\partial^\alpha\Phi_i}(f_i-g_i)\right\|_\A
\leq
\sum_{i=1}^m\|\partial^\alpha\Phi_i\|_{L^1(\Xi)}\| f_i-g_i\|_\A \le\epsilon .
\end{gather*}
Thus we only need to prove that for each $\Phi\in C_c^\infty(\Xi)$ and $g\in\Q(I)\cdot\A$ the element
$\Th_\Phi(g)$ belongs to $\Q(I)\cdot\A^\infty$. Let $\varphi_1,\dots,\varphi_j\in I$ and
$h_1,\dots,h_j\in\A$ such that $g=\sum\limits_{i=1}^j\Q(\varphi_i)h_i$. Then
\[
\Th_\Phi(g)=\sum_{i=1}^j\Th_\Phi [\Q(\varphi_i)h_i ],
\]
and by covariance, for each index $i$ one has
\begin{gather*}
\Th_\Phi\left[\Q(\varphi_i)h_i\right]=\int_\Xi dY\Phi(Y)\Q(\varphi_i)\Th_X(h_i)=\Q(\varphi_i)
\left[\Th_\Phi(h_i)\right]\in\Q(I)\cdot\A^\infty .\tag*{\qed}
\end{gather*}
\renewcommand{\qed}{}
\end{proof}

\begin{Theorem}\label{crifel}
Rieffel quantization transforms covariant $\C(T)$-algebras into covariant $\C(T)$-algebras: there exists a nondegenerate
monomorphism $\mathfrak Q:\C(T)\rightarrow\mathcal Z\mathcal M(\AA)$ satisfying for all $\varphi\in\C(T)$, $f\in\A$ and
$X\in\Xi$ the covariance relation
$\Th_X[\QQ(\varphi)f]=\QQ(\varphi) [\Th_X(f) ]$.
\end{Theorem}

\begin{proof}
The action $\Th$ of $\Xi$ on $\A$ extends canonically to an action by automorphisms of the multiplier algebra
$\mathcal M(\A)$, also denoted by $\Th$, which is not strongly continuous in general.
But, tautologically, it restricts to a strongly continuous action
$\Th:\Xi\rightarrow{\rm Aut}[\mathcal M_0(\A)]$ on the $C^*$-subalgebra
\begin{gather*}
\mathcal M_0(\A):=\{m\in\mathcal M(\A)\mid\Xi\ni X\mapsto\Th_X(m)\in\mathcal M(\A)\ {\rm is\ norm\ continuous}\}.
\end{gather*}
In these terms, the covariance condition on $\Q$ says simply that for any $\varphi\in\C(T)$ the multip\-lier~$\Q(\varphi)$ is a f\/ixed point for all the automorphisms~$\Th_X$
(take $f=1$ in~(\ref{incep})). As a very weak consequence one has
$\Q[\C(T)]\subset\mathcal M_0(\A)^\infty$, with an obvious notation for the smooth vectors.

Proposition~5.10 from \cite{Rie1} applied to the unital $C^*$-algebra $\mathcal M_0(\A)$ says that the Rief\/fel quantization of
$\mathcal M_0(\A)$ is a $C^*$-subalgebra of $\mathcal M(\AA)$. Consequently one has $\Q[\C(T)]\subset\mathcal M_0(\A)^\infty
\subset\mathcal M(\AA)$ and this supplies a candidate $\QQ:\C(T)\rightarrow\mathcal M(\AA)$.
This is obviously an injective map and the range is only composed of f\/ixed points, which insures covariance.

Let us set for a moment $\mathcal M:=\mathcal M_0(\A)$, with multiplication~$\cdot$, and denote by
$\mathfrak M\subset\mathcal M(\AA)$ its Rief\/fel quantization, with multiplication legitimately denoted by~$\#$.
For smooth elements $m,n\in\mathcal M^\infty=\mathfrak M^\infty$, one of them being a f\/ixed point central in~$\mathcal M$,
one has $m\#n=m\cdot n=n\cdot m=n\#m$ (Corollary~2.13 in~\cite{Rie1}).
Thus the mapping~$\QQ$ is again a monomorphism and its range is contained in $\mathcal Z\mathfrak M$.
A density argument with respect to the strict topology implies that every $\mathfrak Q(\varphi)$
commutes with all the elements of~$\mathcal M(\AA)$, thus $\QQ[\C(T)]\subset\mathcal Z\mathcal M(\AA)$as required.

Now we only need to show nondegeneracy, i.e.\ the fact that $\QQ[\C(T)]\cdot\AA$ is dense in~$\AA$.
We show the even stronger assertion that $\Q[\C(T)]\cdot\A^\infty=\QQ[\C(T)]\cdot\AA^\infty$ is dense in~$\AA$.
This would follow if we knew that $\Q[\C(T)]\cdot\A^\infty$ is dense in $\AA^\infty$ with respect to its Fr\'echet topology
given by the semi-norms~(\ref{semon});
then we use denseness of $\AA^\infty$ in the weaker $C^*$-norm topology of~$\AA$.

We recall from Section~\ref{sectra} that $\A^\infty$ and $\AA^\infty$ coincide even as Fr\'echet spaces.
Therefore one is reduced to showing that $\Q[\C(T)]\cdot\A^\infty$ is dense in $\A^\infty$ for its Fr\'echet topology.
Taking $\mathcal I=\C(T)$ in Lemma~\ref{densly}, we f\/ind out that $\Q[\C(T)]\cdot\A^\infty$ is dense in
$\big(\overline{\Q[\C(T)]\cdot\A}\big)\cap\A^\infty$, which equals~$\A^\infty$ since $\Q$ has been assumed
nondegenerate. This f\/inishes the proof.
\end{proof}

If $\A$ is a covariant $\C(T)$-algebra, then $\mathcal I(t):=\overline{\Q\left[\C_t(T)\right]\cdot\A}$
is an invariant ideal of~$\A$. We can apply Rief\/fel quantization to
$\mathcal I(t)$, to $\A(t):=\A/\mathcal I(t)$ (with the obvious actions of $\Xi$) and to the projection
$\mathcal P(t):\A\rightarrow\A(t)$. One gets $C^*$-algebras~$\mathfrak I_t$, $\AA_t$ as well as the
morphism $\mathfrak P_t:\AA\rightarrow\AA_t$. By Theorem~7.7 at~\cite{Rie1} the kernel of~$\mathfrak P_t$ is~$\mathfrak I_t$, so~$\AA_t$ can be identif\/ied to the quotient~$\AA/\mathfrak I_t$.

On the other hand, by using the $\C(T)$-structure of the $C^*$-algebra $\AA$ given by Theorem~\ref{crifel}, we have ideals
$\mathfrak I(t):=\overline{\QQ\left[\C_t(T)\right]\cdot\AA}$ in $\AA$ as well as quotients $\AA(t):=\AA/\mathfrak I(t)$
to which we associates projections $\AA\overset{\PP(t)}{\longrightarrow}\AA(t)$.
However, one gets

\begin{Proposition}\label{pieligroasa}
With notation as above, for each $t\in T$ we have $\mathfrak I(t)=\mathfrak I_t$.
In particular, the fibers $\AA(t)=\AA/\mathfrak I(t)$ of  the $\C(T)$-algebra $\AA$ are isomorphic
to the Rieffel quantization $\AA_t$ of the fibers $\A(t)=\A/\mathcal I(t)$ of $\A$ and for each
$f\in\AA$ the mapping $t\mapsto\| \mathfrak P(t)f\|_{\AA(t)}=\| \PP_t f\|_{\AA_t}$
is upper-semi-continuous.
\end{Proposition}

\begin{proof}
We recall that $\mathcal I(t)^\infty$ and $\mathfrak I(t)^\infty$ coincide as Fr\'echet spaces.
By Lemma \ref{densly}, $\QQ\left[\C_t(T)\right]\cdot\AA^\infty$ is dense in $\mathfrak I(t)^\infty$,
thus in $\mathfrak I(t)$, and
$\Q\left[\C_t(T)\right]\cdot\A^\infty$ is dense in $\mathcal I(t)^\infty=\mathfrak I(t)^\infty$,
thus also dense in~$\mathfrak I_t$.

By construction one has $\QQ\left[\C_t(T)\right]\cdot\AA^\infty=\Q\left[\C_t(T)\right]\cdot\A^\infty$;
consequently $\mathfrak I(t)=\mathfrak I_t$ for every $t\in T$ and the proof is f\/inished.
\end{proof}

\begin{Remark}\label{spekul}
We are going to say that
$\big\{\mathcal A\overset{\P(t)}{\longrightarrow}\mathcal A(t)\mid t\in T\big\}$ and
$\big\{\AA\overset{\PP(t)}{\longrightarrow}\AA(t)\mid t\in T\big\}$ are {\it covariant upper-semi-continuous fields of
$C^*$-algebras}.
The intrinsic def\/inition, in the f\/irst case for instance, would be the following:
$\big\{\mathcal A\overset{\P(t)}{\longrightarrow}\mathcal A(t)\mid t\in T\big\}$ is required to be an upper semi-continuous
f\/ield of $C^*$-algebras and we also ask the action $\Th$ to leave invariant all the ideals
$\mathcal I(t)=\ker[\P(t)]$. It is easily seen that this is equivalent to requiring the covariance of the associated
$\C(T)$-structure. This makes the connection with Def\/inition~3.1 in~\cite{Rie}.
\end{Remark}

For section $C^*$-algebras of an upper-semi-continuous f\/ield it is known \cite{Wi} that each irreducible representation
factorizes through one of the f\/ibers. Therefore we get

\begin{Corollary}\label{factor}
Let $(\A,\Th,\Xi,[\![\cdot,\cdot]\!])$ be a classical data and assume that $\A$ is a
$\Th$-covariant $\C(T)$-algebra with respect to a
Hausdorff locally compact space $T$, with fibers $\{\A(t)\mid t\in T\}$. Denote, respectively, by $\AA$ and $\AA(t)$
the corresponding quantized $C^*$-algebras. Then any irreducible representation of $\AA$ factorizes through
one of the algebras $\AA(t)$.
\end{Corollary}

The $\C(T)$-structure $\QQ$ of $\AA$, given by Theorem \ref{crifel}, def\/ines canonically the map
$\mathfrak q{:}\!$ ${\rm Prim}(\AA){\rightarrow} T$, as explained at the end of Section~\ref{secintra}.
If $\pi:\AA\rightarrow\mathbb B(\mathcal H)$ is an irreducible Hilbert space representation of $\AA$, then
the point $t$ in Corollary \ref{factor} is $\mathfrak q[\ker(\pi)]$.

\section{Lower-semi-continuity under Rief\/fel quantization}\label{joasa}

We keep the previous setting and inquire now if lower-semi-continuity of the mappings
$t\mapsto\| \P(t)f\|_{\A(t)}$ for all $f\in\A$ implies lower-semi-continuity of the mappings
$t\mapsto\| \PP(t)f\|_{\AA(t)}$ for all $f\in\AA$. We start by noticing that
${\rm Prim}(\A)$ and ${\rm Prim}(\AA)$
are canonically endowed with continuous actions of the group $\Xi$; once again these actions will be denoted by~$\Th$.
By the discussion at the end of Section~\ref{secintra} we are left with proving

\begin{Proposition}\label{segund}
Suppose that $\Q:\C(T)\to \mathcal Z\mathcal M(\A)$ is a covariant $\C(T)$-algebra structure on~$\A$ and that the associated function
$q:{\rm Prim}(\A)\to T$ is open. Then the function $\mathfrak q:{\rm Prim}(\AA)\to T$ associated to
$\mathfrak Q:\C(T)\to \mathcal Z(\AA)$ is also open.
\end{Proposition}

\begin{proof}
We remark f\/irst that $q$ is $\Th$-covariant (Lemma~8.1 in~\cite{Wi}), i.e.\ one has $q\circ\Th_X=q$ for every $X\in\Xi$.
Consequently, if $\mathcal O\subset {\rm Prim}(\A)$ is
an open set, then $\Th_\Xi(\mathcal O):=\{\Th_X(\mathcal K)\mid X\in\Xi, \, \mathcal K\in\mathcal O\}$ will also be an open set
and $q(\mathcal O)=q\left[\Th_\Xi(\mathcal O)\right)]$.
So $q$ will be open if\/f it sends open {\it invariant} subsets of ${\rm Prim}(\A)$
into open subsets of $T$. The same is true for $\mathfrak q:{\rm Prim}(\AA)\rightarrow T$. But the most general
open subset of ${\rm Prim}(\A)$ has the form
\[
\mathcal O_{\mathcal J}:=\{\mathcal K\in {\rm Prim}(\A)\mid \mathcal J\not\subset\mathcal K\}=h(\mathcal J)^c
\]
for some ideal $\mathcal J$ of $\A$, being the complement of the hull  $h(\mathcal J)$ of this ideal.
In addition, $\mathcal O_{\mathcal J}$ is $\Th$-invariant if\/f~$\mathcal J$ is an
invariant ideal. We also recall that Rief\/fel quantization establishes
a~one-to-one correspondence between invariant ideals of $\A$ and invariant ideals of~$\AA$.

So let $\mathcal J$ be an invariant ideal in $\A$ and $\mathfrak J$ its quantization
(an invariant ideal in~$\AA$). We would like to show that
$q\left(\mathcal O_\mathcal J\right)=\mathfrak q\left(\mathcal O_\mathfrak J\right)$;
by the discussion above this would imply that $q$ and $\mathfrak q$ are simultaneously open. Using the fact that
$q(\mathcal K)=t$ if and only if $\mathcal I(t)\subseteq \mathcal K$ and similarly for $\mathfrak q$, one gets
\[
q\left(\mathcal O_\mathcal J\right)=\{t\in T\mid \exists\,\mathcal K\in {\rm Prim}(\A),\ \mathcal J\not\subset\mathcal K,
\ \mathcal I(t)\subset\mathcal K\}
\]
and
\[
\mathfrak q\left(\mathcal O_\mathfrak J\right)=\{t\in T\mid \exists\,\mathfrak K\in {\rm Prim}(\AA),\
\mathfrak J\not\subset\mathfrak K,\ \mathfrak I(t)\subset\mathfrak K\}.
\]
Using the hull application  and the fact that both the hull and the kernel are decreasing, one can write
\[
t\notin q\left(\mathcal O_\mathcal J\right)\ \Longleftrightarrow\ h[\mathcal I(t)]\cap h[\mathcal J]^c=\varnothing
\ \Longleftrightarrow\ h[\mathcal I(t)]\subset h[\mathcal J]\ \Longleftrightarrow\ \mathcal I(t)\supset\mathcal J
\]
and
\[
t\notin \mathfrak q\left(\mathcal O_\mathfrak J\right)\ \Longleftrightarrow\ h[\mathfrak I(t)]\cap h[\mathfrak J]^c=\varnothing
\ \Longleftrightarrow\ h[\mathfrak I(t)]\subset h[\mathfrak J]\ \Longleftrightarrow\ \mathfrak I(t)\supset\mathfrak J.
\]
To f\/inish the proof one only needs to notice that the Rief\/fel quantization of invariant ideals preserves inclusions.
\end{proof}

\begin{Remark}\label{harto}
The def\/inition of {\it a covariant continuous field of $C^*$-algebras} is naturally obtained by adding the lower-semi-continuity condition to the def\/inition of an upper-semi-continuous f\/ield of $C^*$-algebras contained in Remark~\ref{spekul}. Using this notion, {\it Theorem~{\rm \ref{mein}} is now fully justified.}
Another proof, based on the theory of twisted crossed products and on results from the deep work~\cite{Rie} can be found in~\cite{BM}, as we indicated in the introduction.
\end{Remark}

The $C^*$-dynamical system $(\A,\Th,\Xi)$ being given, one could try one of the choices $T={\rm Orb}[{\rm Prim}(\A)]$
({\it the orbit space}) or $T={\rm Quorb}[{\rm Prim}(\A)]$ ({\it the quasi-orbit space}), both associated to the natural
action of $\Xi$ on the space ${\rm Prim}(\A)$. We recall that, by def\/inition, {\it a quasi-orbit} is the closure of an orbit
and we refer to \cite{Wi} for all the fairly standard assertions we are \mbox{going} to make about these spaces.
The two spaces are quotients of ${\rm Prim}(\A)$ with respect to obvious equivalence relations.
Endowed with the quotient topology they are locally compact, but they may fail to possess the Hausdorf\/f property.
But {\it the orbit map} $p:{\rm Prim}(\A)\rightarrow{\rm Orb}[{\rm Prim}(\A)]$ and
{\it the quasi-orbit map} $q:{\rm Prim}(\A)\rightarrow{\rm Quorb}[{\rm Prim}(\A)]$ are continuous open surjections.
So one can state:

\begin{Corollary}\label{nech}
If the quasi-orbit space associated to the dynamical system
$({\rm Prim}(\A),\Th,\Xi)$ is Hausdorff, then the deformed $C^*$-algebra $\AA$ can be expressed as a continuous field of
$C^*$-algebras over the base ${\rm Quorb}[{\rm Prim}(\A)]$.
A similar statement holds with ``quasi-orbit'' replaced by ``orbit'' and ${\rm Quorb}[{\rm Prim}(\A)]$ replaced
by ${\rm Orb}[{\rm Prim}(\A)]$.
\end{Corollary}

Notice that, when ${\rm Orb}[{\rm Prim}(\A)]$ happens to be Hausdorf\/f, the orbits will be automatically closed
(as inverse images by $p$ of points); so
one would actually have ${\rm Orb}[{\rm Prim}(\A)]={\rm Quorb}[{\rm Prim}(\A)]$.

\section{Some examples}\label{seretide}

The most important is the Abelian case, that has been introduced at the end of Section \ref{sectra}.
We make a brief review of this case in conjunction to continuity matters, following the more detailed version of \cite[Section~2]{BM}. This will illustrate the general theory and will also be a preparation for some of the examples we intend to present below.

We recall that an action $\Th$ of $\Xi$ by homeomorphisms of the locally compact space $\Si$ is given. We also assume given a continuous surjection $q:\Sigma\to T$. Then we have the disjoint decomposition of $\Si$ in closed subsets
\begin{gather*}
\Si=\sqcup_{t\in T}\Si_t,\qquad \Si_t:=q^{-1}(\{t\}).
\end{gather*}
Associated to the canonical injections $j_t:\Si_t\to\Si$, we have associated restriction epimorphisms
\begin{gather*}
\R(t): \ \C(\Si)\to\C(\Si_t),\qquad \R(t)f:=f|_{\Si_t}=f\circ j_t,\quad \forall\, t\in T.
\end{gather*}
We say that the continuous surjection $q$ {\it is $\Th$-covariant} if each $\Si_t$ is $\Theta$-invariant.

The Rief\/fel-quantized $C^*$-algebras $\CC(\Si)$ and $\CC(\Si_t)$ as well as the epimorphisms
$\mathfrak R(t):\CC(\Si)\rightarrow\CC(\Si_t)$ were introduced above.
Applying now the results obtained in Sections~\ref{isra} and~\ref{joasa}, one gets as in \cite[Section~2]{BM}

\begin{Corollary}\label{cfild}
Assume that the mapping $q:\Si\rightarrow T$ is a $\Th$-covariant continuous surjection.
Then $\big\{\CC(\Si)\overset{\RR(t)}{\longrightarrow}\CC(\Si_t)\mid  t\in T\big\}$ is a covariant
upper-semi-continuous field of noncommutative $C^*$-algebras.
If~$q$ is also open, then the field is continuous.
\end{Corollary}

Let us assume now that the orbit space ${\rm Orb}(\Si)$ is Hausdorf\/f.
Any orbit, being the inverse image of a point in ${\rm Orb}(\Si)$, will be closed in
$\Si$ and invariant; it will also be homeomorphic to the quotient of $\Xi$ by the corresponding
stability group. As a precise particular case of Corollary~\ref{nech} one can state:

\begin{Corollary}\label{neluz}
If the orbit space of the dynamical system
$(\Si,\Th,\Xi)$ is Hausdorff, then the deformed $C^*$-algebra $\CC(\Si)$ can be expressed as a continuous field of
$C^*$-algebras over the base space ${\rm Orb}(\Si)$. The fiber over $\mathcal O\in{\rm Orb}(\Si)$ is the deformation
of the Abelian algebra $\C(\mathcal O)\cong\C(\Xi/\Xi_\mathcal O)$.
\end{Corollary}

\begin{Remark}\label{marcu}
It is known that the orbit space is Hausdorf\/f if the action $\Th$ is proper; this happens for instance if~$\Si$ is a Hausdorf\/f locally compact group on which the closed subgroup $\Xi$ acts by
left translations. More generally, assume that the action $\Th$ factorizes
through a compact group~$\widehat \Xi$,
i.e.\ the kernel of~$\Th$ contains a closed co-compact subgroup~$Z$ of $\Xi$ (with $\widehat \Xi=\Xi/Z$).
Then the orbit space under the initial action is the same as the orbit space of the
action of the compact quotient. But the action of a compact group is proper and Corollary~\ref{neluz} applies.
\end{Remark}

\begin{Example}\label{primu}
Let $\mathscr A$ be a $C^*$-algebra and $T$ a locally compact space. On
\begin{gather*}
\A\equiv\C(T;\mathscr A):=\{f:T\rightarrow\mathscr A\mid f\ \text{is continuous and small at inf\/inity}\}
\end{gather*}
we consider the natural structure of $C^*$-algebra. It clearly def\/ines a continuous f\/ield of $C^*$-algebras
\[
\big\{\C(T;\mathscr A)\overset{\delta(t)}{\longrightarrow}\mathscr A\mid t\in T\big\},\qquad \delta(t)f:=f(t).
\]
The associated $\C(T)$-structure is given by $\left[\Q(\varphi)f\right](t):=\varphi(t)f(t)$ for $\varphi\in\C(T)$, $f\in\A$, $t\in T$.
For each $t\in T$ an action $\th^t$ of $\Xi$ on $\mathscr A$ is given; we require that each $f\in\A$ verif\/ies
\begin{gather*}
\sup_{t\in T}\| \th^t_X[f(t)]-f(t) \|_\mathscr A\underset{X\rightarrow 0}{\longrightarrow} 0.
\end{gather*}
Then obviously
\begin{gather*}
\Th: \ \Xi\rightarrow{\rm Aut}(\A),\qquad [\Th_X(f)](t):=\th^t_X[f(t)]
\end{gather*}
def\/ines a continuous action of the vector group $\Xi$ on $\A$. Each of the kernels
\[
\mathcal I(t):=\ker[\delta(t)]=\{f\in\C(T;\mathscr A)\mid f(t)=0\}
\]
is $\Th$-invariant, so one actually has a covariant continuous f\/ield of $C^*$-algebras
(see Remarks~\ref{spekul} and~\ref{harto}).
It makes sense to apply Rief\/fel quantization, getting $C^*$-algebras (respectively)
$\AA\equiv\CC(T;\mathscr A)$ from the dynamical system $(\A\equiv\C(T;\mathscr A),\Th)$
and $\AA(t)$ from the dynamical system $(\mathscr A,\th^t)$ for all $t\in T$. From the results above one concludes that
$\big\{\AA\overset{\Delta(t)}{\longrightarrow}\AA(t)\mid t\in T\big\}$ {\it is also a covariant field of $C^*$-algebras}.
For each $t$ we denoted by $\Delta(t)$ the Rief\/fel quantization of the morphism~$\delta(t)$.
\end{Example}

\begin{Example}\label{secundu}
A particular case, considered in \cite[Chapter~8]{Rie1}, consists in taking $T:={\rm End}(\Xi)$ the space of all linear maps
$t:\Xi\rightarrow\Xi$; it is a locally compact (f\/inite-dimensional vector) space with the obvious operator norm.
If an initial action~$\theta$ of~$\Xi$ on $\mathscr A$ is f\/ixed, the choice $\theta^t_X:=\theta_{tX}$ verify all
the requirements above. Therefore one gets a covariant continuous f\/ield of $C^*$-algebras indexed by ${\rm End}(\Xi)$.
This is basically \cite[Theorem~8.3]{Rie1}; we think that our treatment gives a~simpler and more unif\/ied proof of this result,
especially concerning the lower-semi-continuous part.
In particular, for any $f\in\C[{\rm End}(\Xi);\mathscr A]$, one has
$\lim\limits_{t\rightarrow 0}\| f(t)\|_{\AA(t)}=\| f(0)\|_\mathscr A$.
An interesting particular case is obtained restricting the arguments to the compact subspace
$T_0:=\{t=\sqrt\hb\, {\rm id}_\Xi\mid \hb\in[0,1]\}\subset T$. The number $\hb$ corresponds to the Plank constant and, even for
constant $f:[0,1]\rightarrow\mathscr A$, the relation
$\lim\limits_{\hb\rightarrow 0} \| f\|_{\AA(\hb)}=\| f\|_\mathscr A$
is nontrivial and has an important physical
interpretation concerning the semiclassical behavior of the quantum mechanical formalism. We refer to \cite{La,Rie1,Rie2} for
much more on this topic.
\end{Example}

\begin{Remark}\label{subsumat}
A way to convert Example \ref{primu} into a more sophisticated one is as follows:

For every $t\in T$ pick $\B(t)$ to be a $C^*$-subalgebra
of $\mathscr A$ which is invariant under the action~$\th^t$. Construct the $C^*$-subalgebra~$\B$ of~$\A$ def\/ined as
$\B:=\{f\in\C(T;\mathscr A)\mid f(t)\in\B(t),\ \forall\,t\in T\}$,
which is obviously invariant under the action $\Th$. One gets a covariant continuous f\/ield of $C^*$-algebras
$\big\{\B\overset{\P(t)}{\longrightarrow}\B(t)\mid t\in T\big\}$, where $\P(t)$ is a restriction of the
epimorphism $\delta(t)$. The general theory developed in Sections~\ref{isra} and~\ref{joasa}
supplies another covariant continuous f\/ield of $C^*$-algebras
$\big\{\BB\overset{\PP(t)}{\longrightarrow}\BB(t)\mid t\in T\big\}$, where~$\BB(t)$ is the quantization of~$\B(t)$
and can be identif\/ied with an invariant $C^*$-subalgebra of~$\AA(t)$.
\end{Remark}

\begin{Example}\label{hartog}
Crossed products associated to actions of $\X:=\mathbb R^n$ on $C^*$-algebras can be obtained from Rief\/fel's
quantization procedure, as it is explained in \cite[Example~10.5]{Rie1}. From the results of the present article one could
infer rather easily, as a particular case, that (informally) {\it the crossed product by a continuous field of
$C^*$-algebras is a continuous field of crossed products}. Such results exist in a much greater generality, including
(twisted) actions of amenable locally compact groups~\cite{Ni,PR,Rie,Wi}, so we are not going to give details.
\end{Example}

\begin{Example}\label{tertu}
In \cite{Rie3} one constructs $C^*$-algebras which can be considered quantum versions of a certain class of compact
connected Lie groups. We will have nothing to say about the extra structure making them quantum groups; we are only going to
apply the results above to present these $C^*$-algebras as continuous f\/ields.

Let $\Si$ be a compact connected Lie group, containing {\it a toral subgroup}, i.e.\ a connected closed Abelian subgroup $H$.
Such a toral group is isomorphic to an $n$-dimensional torus $\mathbb T^n$. Assume given a continuous group epimorphism
$\eta:\mathbb R^n\rightarrow H$ (for example the exponential map def\/ined on the Lie algebra $\mathfrak H\equiv\mathbb R^n$).
We use $\eta$ to def\/ine an action of $\Xi:=\mathbb R^n\times\mathbb R^n$ on $\Si$ by $\Th_{(x,y)}(\si):=\eta(-x)\si\eta(y)$.
Then, by applying Rief\/fel deformation to $\A:=\C(\Si)$ using the action $\Th$ (and a certain type of skew-symmetric
operator on $\Xi$), one gets the $C^*$-algebra $\AA:=\CC(\Si)$ which, endowed with suitable extra structure, is
regarded as a quantum group corresponding to $\Si$.

It is obvious that the action factorizes through the compact group $H\times H$.
Thus the orbit space ${\rm Orb}(\Si)$ is Hausdorf\/f and
Remark \ref{marcu} and Corollary \ref{neluz} serve to express $\CC(\Si)$ as a continuous f\/ield of $C^*$-algebras.
For the stability group of any orbit $\mathcal O$ one can write
$\Xi_\mathcal O\supset\ker(\Th)\supset \ker(\eta)\times\ker(\eta)$,
thus $\mathcal O\cong\Xi/\Xi_\mathcal O$ is a continuous image of $H\times H$.
\end{Example}

\begin{Example}\label{quartu}
An interesting particular case, taken from \cite{Rie3}, involves the construction of a quantum version of the compact Lie group
$\Si:=\mathbb T\times SU(2)$. Here $\mathbb T$ is the $1$-torus, the group $SU(2)$ contains diagonally a second
copy of $\mathbb T$ and can be parametrised
by the $3$-sphere $S^3:=\{(z,w)\in \mathbb C^2\mid z^2+w^2=1\}$, and so $\Si$ contains a $2$-torus.
Initially $\Xi=\mathbb R^4$ acts on $\Si$ in the
given way, but it is shown in \cite{Rie3} (using results from \cite{Rie1}) that the same deformed algebra is obtained
by the action
\[
\Th': \ \Xi':=\mathbb R^2\rightarrow{\rm Homeo}\left(\mathbb T\times S^3\right),\qquad
\Th'_{(x,y)}(\eta;z,w):=\left(e^{-2\pi ix}\eta;z,e^{4\pi iy}w\right).
\]
The orbit space is homeomorphic with the closed unit disk $T:=\{z\in\mathbb C\mid |z|\le 1\}$.
The orbits corresponding to $|z|<1$ are $2$-tori, while the orbits corresponding to $|z|=1$ (implying $w=0$) are $1$-tori.
If we set $\A:=\C(\mathbb T\times SU(2))$, then the quantized $C^*$-algebra $\AA\cong\CC(\mathbb T\times SU(2))$
deserves to be called {\it a quantum} $\mathbb T\times SU(2)$. The deformation of the continuous functions on
any of the $2$-tori leads to a quantum torus. By multiplying the initial skew-symmetric
form $[\![\cdot,\cdot]\!]$ with an irrational
number~$\beta$ one can make this noncommutative torus $\CC_\beta(\mathbb T^2)$ irrational,
which serves to show that the corresponding quantum
$\mathbb T\times SU(2)$ (obtained for such a $\beta$) is not of type~I. But applying the results obtained here one also gets the
detailed information: {\it The algebra $\CC(\mathbb T\times SU(2))$ can be written over the closed unit disk $T$ as
a continuous field of noncommutative $2$-tori and Abelian $C^*$-algebras $($corresponding to the one-dimensional orbits$)$.}

Many other particular cases can be worked out in detail. We propose to the reader the example $\Si:=SU(2)\times SU(2)$.
\end{Example}

\subsection*{Acknowledgements}

The authors were partially supported by \textit{N\'ucleo
Cientifico ICM P07-027-F ``Mathematical Theo\-ry of Quantum and
Classical Magnetic Systems''}.

\pdfbookmark[1]{References}{ref}
\LastPageEnding

\end{document}